\documentclass[12pt]{my_elsarticle}

\usepackage{amsmath,amssymb,amsthm,amsfonts}

\newcommand{\R}{\mathbb{R}}
\newcommand{\N}{\mathbb{N}}

\newtheorem{theorem}{Theorem}[section]

\newtheorem{example}{Example}[section]

\newtheorem{lemma}[theorem]{Lemma}
\newtheorem{proposition}{Proposition}[section]

\newtheorem{definition}[theorem]{Definition}
\newtheorem{remark}{Remark}[section]

\begin{document}

\begin{frontmatter}

\title{Modified Optimal Energy and Initial Memory\\
of Fractional Continuous-Time Linear Systems\footnote{Submitted 30/Nov/2009; 
Revised (major revision) 24/April/2010; Accepted (after minor revision) 
22/July/2010; for publication in \emph{Signal Processing}.}}

\author[labelDorota]{Dorota Mozyrska}
\ead{d.mozyrska@pb.edu.pl}

\author[labelDelfim]{Delfim F. M. Torres}
\ead{delfim@ua.pt}

\address[labelDorota]{Bia{\l}ystok  University of Technology,
    ul. Wiejska 45A,
    15-351 Bia{\l}ystok, Poland}

\address[labelDelfim]{Department of Mathematics,
   University of Aveiro,
   3810-193 Aveiro, Portugal}

\begin{abstract}
Fractional systems with Riemann-Liouville derivatives
are considered. The initial memory value problem is posed and studied.
We obtain explicit steering laws with respect to the values
of the fractional integrals of the state variables.
The Gramian is generalized and steering functions
between memory values are characterized.
\end{abstract}


\end{frontmatter}


\section{Introduction}

Several approaches to the generalization of the notion of differentiation
to non-integer orders are available in the literature of fractional calculus.
Those include, \textrm{e.g.}, Riemann-Liouville, Gr\"{u}nwald-Letnikow,
Caputo, and the generalized functions approach.
Mathematical fundamentals of fractional calculus are given in the monographs
\cite{miller,nishimoto,oldham,oustalup,podlubny,samko}.
The idea of fractional calculus is used by engineers
to model various processes since late 1960's
\cite{hilfer,sabatier1,scalas,vinagre1,vinagre2}.

Here we study fractional-order systems using a vector
state space representation, where we include also initialization.
Our systems are described by initialized Riemann-Liouville derivatives
\cite{hartley2,hartley1}. Fractional state space representations of systems are given,
\textrm{e.g.}, in \cite{bagley,bettayeb,hartley1,kaczorek1,sierociuk}.
Particularly in \cite{hartley1}, the method of initialized fractional calculus
is addressed \cite{hartley0,hartley,hartley1,lorenzo}.

The initial condition problem for fractional linear systems is a subject under strong consideration
\cite{lorenzo0,lorenzo,lorenzo1,lorenzo2,ortgueira_signal1,ortgueira_signal2,ortigueira2,sabatier}.
Several authors claim that the Riemann-Liouville derivative leads to initial conditions without physical meaning.
Others, \textrm{e.g.} Heymans and Podlubny \cite{heymans_podlubny},
contradict such claim, giving several examples of physical meaning to the Riemann-Liouville
initial conditions, by introducing the concept of ``inseparable twins''.
Proper initialization is crucial to the solution and understanding
of fractional differential equations. An important tool in the solution
of fractional differential or differintegral equations is provided
by the Laplace transform. We consider the problem of solving initialized fractional
differential systems using the Laplace transform.
This can be interpreted as the initial memory of the system.
Then we develop the idea by looking at a given system not
in the perspective of the state but also throughout memory,
that can be of a different fractional order than the order of the system.

The paper is organized as follows. Section~\ref{sec:prel} gives
a brief review of fractional integrals, the fractional Riemann-Liouville derivative
of scalar functions, properties for special functions in fractional calculus, as the
Mittag-Leffler function and its generalizations, the matrix case,
and the $\alpha$-exponential matrix. Then, in Section~\ref{sec:syst:Gram},
we formulate the fractional system $\Sigma_{\alpha,\psi}$ under consideration,
where $\alpha$ denotes the fractional order of the system and $\psi$ the initialization
function. We prove the form of solution for such kind of systems, we introduce the notion
of $\beta$-controllability, and we give a generalization and description
of the controllability Gramian. Our notion of $\beta$-controllability
is strongly connected with the system and the given
initialization function $\psi$.
The new possibilities are, essentially, present
in the special case when $1-\alpha = \beta$. Indeed,
for $1-\alpha = \beta = 0$ we obtain the classical
linear dynamics but we are comparing the integral of the state.
Moreover, for integer-order systems ($\alpha=1$) the initialization term
in the Laplace transform is the Laplace transform of a constant and becomes
$\Psi(s)=a/s$; for fractional-order systems this term is time-varying
into the future, that is, the history of $x(t)$ has the interference
into the infinite future \cite{hartley1}.
The classical notion of controllability \cite{zabczyk}
is obtained from our definition by taking $\alpha=1$ and $\beta=0$,
what is interpreted as controllability without memory.
Our notion is, however, more than a generalization
to fractional-order systems: it introduces a new meaning,
where we do not observe the state of the system
but the integral of state along time.
Finally, in Section~\ref{sec:rank} we state
the notion of steering control and prove its optimality.
We end with Section~\ref{sec:cnc} of
conclusions and some future perspectives.


\section{Preliminaries}
\label{sec:prel}

In this section we make a review of notations, definitions, and some preliminary
facts which are useful for the paper. Since several different notations
and various definitions are used for fractional derivatives, our goal here
is to make precise the definitions and formulas for solutions of linear differential
fractional systems with Riemann-Liouville derivatives and initial conditions
given by the value of a fractional integral at $t_0$.
We recall  definitions of fractional integrals of arbitrary order,
the Riemann-Liouville derivative of order $\alpha\in (0,1)$,
and a description of special functions in the fractional calculus.
We present also the idea of initialization response with connection
to the initialized linear fractional differential equations for order
$\alpha\in (0,1)$. As in the classical (non-fractional) linear case,
the Laplace transform provides a strong tool. We give here
some useful formulas connected with the matter.

\begin{definition}[\cite{kilbas,podlubny,samko}]
Let $\varphi\in L_1\left([t_0,t_1],\R\right)$. The integrals
\begin{equation*}
\label{eq:int1}
I^{\alpha}_{t_0+}\varphi(t)=\frac{1}{\Gamma(\alpha)}
\int_{t_0}^t \varphi(\tau)(t-\tau)^{\alpha-1}d\tau\,, \ \ t>t_0,
\end{equation*}
\begin{equation*}
\label{eq:int2}
I^{\alpha}_{t_1-}\varphi(t)=\frac{1}{\Gamma(\alpha)}
\int_{t}^{t_1} \varphi(\tau)(\tau-t) ^{\alpha-1}d\tau\,,  \ \ t< t_1\, ,
\end{equation*}
where $\Gamma$ is the gamma function and $\alpha>0$,
are called, respectively, the \emph{left-sided} and the \emph{right-sided} fractional integrals
of order $\alpha$.  Additionally, we define
$I^0_{t_0+}=I^0_{t_1-}:=I\hspace{-1.8mm}I$ (identity operator).
\end{definition}

We have the following formula of integration by parts for fractional integrals.

\begin{proposition}[\cite{kilbas}]
Let $\alpha>0$ and $1/p + 1/q \leq 1+\alpha$, $p\geq 1$, $q\geq 1$,
with $p\neq 1$ and $q\neq 1$ in the case $1/p + 1/q = 1+\alpha$.
For $\varphi \in L_p\left([t_0,t_1],\R\right)$ and $\psi\in L_q\left([t_0,t_1],\R\right)$ it holds:
\begin{equation}
\label{intbypart1}
\int_{t_0}^{t_1}\varphi(\tau)I^{\alpha}_{t_0+}\psi(\tau)d\tau
=\int_{t_0}^{t_1}\psi(\tau)I^{\alpha}_{t_1-}\varphi(\tau)d\tau\,.
\end{equation}
\end{proposition}

\begin{remark}[\cite{bhn}]
When $t_0=0$ we write $I^{\alpha}_{0+}=I^{\alpha}$ and then
$I^{\alpha}f(t)=(f*\varphi_{\alpha})(t)$, where
$\varphi_{\alpha}(t)=\frac{t^{\alpha-1}}{\Gamma(\alpha)}$ for $t>0$,
$\varphi_{\alpha}(t)=0$ for $t\leq 0$,
and $\varphi_{\alpha}\rightarrow \delta(t)$ as $\alpha\rightarrow 0$,
with $\delta$ the delta Dirac pseudo function.
\end{remark}

\begin{definition}[\cite{kilbas,podlubny}]
Let $\varphi$ be defined on the interval $[t_0,t_1]$.
The left-sided Riemann-Liouville derivative of order $\alpha$
and the lower limit $t_0$ are defined through the following:
\begin{equation*}
\label{eq:rl}
D^{\alpha}_{t_0+}\varphi(t)=\frac{1}{\Gamma(n-\alpha)} \left(\frac{d}{dt}\right)^n
\int_{t_0}^t  \varphi(\tau)(t-\tau)^{n-\alpha-1}d\tau \, ,
\end{equation*}
where $n$ is a natural number satisfying $n=[\alpha]+1$
with $[\alpha]$ denoting the integer part of  $\alpha$.
Similarly, the right-sided Riemann-Liouville derivative of order $\alpha$
and the upper limit $t_1$ are defined by
\begin{equation*}
\label{eq:rl2}
D^{\alpha}_{t_1-}\varphi(t)=\frac{1}{\Gamma(n-\alpha)} \left(-\frac{d}{dt}\right)^n
\int_{t}^{t_1}  \varphi(\tau)(\tau-t)^{n-\alpha-1}d\tau \, .
\end{equation*}
\end{definition}

\begin{remark}
Let $\alpha\in(0,1)$. For functions $\varphi$ given
in the interval $[t_0,t_1]$, each of the following expressions
\begin{equation*}
D^{\alpha}_{t_0+}\varphi(t)=\frac{1}{\Gamma(1-\alpha)} \frac{d}{dt}
\int_{t_0}^t  \varphi(\tau)(t-\tau)^{-\alpha}d\tau=  \frac{d}{dt}\left(I^{1-\alpha}_{t_0+}\varphi(t)\right) \, ,
\end{equation*}
\begin{equation*}
D^{\alpha}_{t_1-}\varphi(t)=-\frac{1}{\Gamma(1-\alpha)} \frac{d}{dt}
\int_{t}^{t_1}  \varphi(\tau)(\tau-t)^{-\alpha}d\tau = -  \frac{d}{dt}\left(I^{1-\alpha}_{t_1-}\varphi(t)\right)\,
\end{equation*}
are called  a fractional derivative of order $\alpha$, left-sided and right-sided respectively.
\end{remark}

From \cite[Theorem 2.4]{kilbas} we have the following properties:

\begin{proposition}
If $\alpha>0$, then $D^{\alpha}_{t_0+}I^{\alpha}_{t_0+}\varphi(t)
=\varphi(t)$ for any $\varphi\in L_1(t_0,t_1)$, while
$I^{\alpha}_{t_0+}D^{\alpha}_{t_0+}\varphi(t)=\varphi(t)$
is satisfied for $\varphi\in I^{\alpha}_{t_0+}(L_1(t_0,t_1))$ with
\[I^{\alpha}_{t_0+}(L_1(t_0,t_1))=\{\varphi(t):
\ \varphi(t)=I^{\alpha}_{t_0+}\psi(t), \ \psi \in L_1(t_0,t_1)\} \, .\]
However,
\[I^{\alpha}_{t_0+}D^{\alpha}_{t_0+}\varphi(t)=
\varphi(t)-\sum_{k=0}^{n-1}\frac{(t-t_0)^{\alpha-k-1}}{\Gamma(\alpha
-k)}\left(I^{n-\alpha}_{t_0+}\varphi(t)\right)_{t=t_0}\, .\]
In particular, for $\alpha\in (0,1]$ we have
\[I^{\alpha}_{t_0+}D^{\alpha}_{t_0+}\varphi(t)=\varphi(t)
- \frac{(t-t_0)^{\alpha-1}}{\Gamma(\alpha)}\left(I^{1-\alpha}_{t_0+}\varphi(t)\right)_{t=t_0}\,.\]
\end{proposition}

Similar results hold for right-sided fractional derivatives.
There exist also a formula for the composition of
Riemann-Liouville derivatives (\textrm{cf.} \cite{podlubny}).

The next proposition is based on \cite[Corollary 2, p.~46]{samko}
and is particularly useful for our purposes.

\begin{proposition}
\label{intbypart2}
The formula
\begin{equation}\label{eq:intbypart}
\int_{t_0}^{t_1}f(t)D^{\alpha}_{t_0+}g(t)dt
=\int_{t_0}^{t_1}g(t)D^{\alpha}_{t_1-}f(t)dt\, ,  \ \  0<\alpha<1 \, ,
\end{equation}
is valid under the assumption that $f(t)\in I^{\alpha}_{t_1-}(L_p)$
and $g(t)\in I^{\alpha}_{t_0+}(L_q)$ with $1/p+1/q\leq 1+\alpha$.
\end{proposition}

\begin{proof}
The equality (\ref{eq:intbypart}) follows from (\ref{intbypart1})
if we denote $D^{\alpha}_{t_1-}f(t)=\varphi(t)$, $D^{\alpha}_{t_0+}g(t)=\psi(t)$,
and take into account that $I^{\alpha}_{t_0+}D^{\alpha}_{t_0+}f(t)=f(t)$
is valid for $f(t)\in I^{\alpha}_{t_0+}(L_1)$ (\textrm{cf.} \cite{samko}).
\end{proof}

Another important result is the formula for the Laplace transformation of the derivative
of a function $\varphi$ (\textrm{cf.} \cite{debnath}).
Like the Laplace transformation of an integer order derivative, it is easy
to show that the Laplace transformation of a fractional order
Riemann-Liouville derivative is given by
$$
\mathcal{L}\left[D^{\alpha}_{0+}\varphi(t)\right](s)
=s^{\alpha}\mathcal{L}\left[\varphi(t)\right](s)
- \sum_{k=0}^{n-1}c_ks^k, \  \ n-1< \alpha\leq n\, ,
$$
where $c_k=\left(I^{k+1-\alpha}_{0+}\varphi(t)\right)_{t=0}$.
For $\alpha\in(0,1]$, we have $
\mathcal{L}\left[D^{\alpha}_{0+}\varphi(t)\right](s)=s^{\alpha}\mathcal{L}\left[\varphi(t)\right](s)-a$,
where according with many authors the condition $a=\left(I^{1-\alpha}_{0+}\varphi(t)\right)_{t=0}$
is considered a technical initial condition, mainly due to the Laplace transform
and without a good physical interpretation
\cite{hartley0,hartley,hartley2,hartley1,lorenzo0,lorenzo,lorenzo1,lorenzo2,ortgueira_signal1,ortgueira_signal2,ortigueira2,sabatier}.

In papers of Hartley and Lorenzo \cite{hartley0,hartley,hartley2,hartley1},
one can find the fundamentals about the idea of initialization process.
By $d^{\alpha}_t\varphi$ we denote the derivative without initialization
or just the Liouville derivative defined by the inverse Laplace transform:
\[\mbox{}_0d^{\alpha}_t\varphi(t):=\mathcal{L}^{-1}\left[s^{\alpha}F(s)\right](t)\, ,\]
where $F(s)=\mathcal{L}\left[\varphi(t)\right](s)$.
For a fixed function $\varphi$ and $\alpha>0$,
by initialization function $\psi=\psi(t)$
we mean a function that satisfies the following:
\[D^{\alpha}_{0+}\varphi(t)=\mbox{}_0d_t^{\alpha}(\varphi)(t)+\psi(t)\,.\]
Here we consider only problems starting at 0. Then we have
\[\mathcal{L}\left[D^{\alpha}_{0+}\varphi(t)\right](s)
=s^{\alpha}F(s)+\mathcal{L}\left[\psi(t)\right](s)\,.\]

For a function $x:[0,T]\rightarrow\R^n$, we use similar notation as in the classical case:
\[
D^{\alpha}_{0+}x(t)=D^{\alpha}_{0+}
\left(
\begin{array}{c}
x_1(t)\\ x_2(t)\\ \vdots \\ x_n(t)
\end{array}
\right)
= \left(
\begin{array}{c}
D^{\alpha}_{0+}x_1(t)\\
D^{\alpha}_{0+}x_2(t)\\
\vdots \\
D^{\alpha}_{0+}x_n(t)
\end{array}
\right) \,.
\]
Such situation, when for each component we use the same fractional order $\alpha$ of differentiation,
is called in the literature the fractional-order derivative with commensurate order.
Commensurate order is used
in connection with both Riemann-Liouville
or the Caputo derivative (\textrm{cf.} \cite{bettayeb,kaczorek1,sierociuk}).

We need to mention now two important functions and their extensions to matrices.

\begin{definition}[\cite{eld}]
\label{def:ML:2p}
The two-parameter Mittag-Leffler type function is defined by the series expansion
\begin{equation*}
E_{\alpha,\beta}(z)=\sum_{k=0}^{\infty}\frac{z^{k}}{\Gamma(k\alpha +\beta)}\, , \quad
\alpha>0 \, , \quad \beta>0\,.
\end{equation*}
Let $A\in\R^{n\times n}$. By
\begin{equation*}
E_{\alpha,\beta}(At^{\alpha})=\sum_{k=0}^{\infty}A^k\frac{t^{k\alpha}}{\Gamma(k\alpha +\beta)}
\end{equation*}
it is denoted the extension of the two-parameter Mittag-Leffler function to matrices.
\end{definition}

For $\beta=1$ we obtain the Mittag-Leffler function of one parameter:
\begin{equation*}
E_{\alpha}(z)=\sum_{k=0}^{\infty}\frac{z^{k}}{\Gamma(k\alpha +1)} \, .
\end{equation*}
In the matrical case, $E_{\alpha}(At^{\alpha})
=\sum_{k=0}^{\infty}A^k\frac{t^{k\alpha}}{\Gamma(k\alpha +1)}$.

\begin{definition}[\cite{kilbas}]
\label{def:aEM}
Let $A\in \R^{n\times n}$. Then,
\begin{equation*}
e_{\alpha}^{At}=t^{\alpha-1}\sum_{k=0}^{\infty}A^k\frac{t^{k\alpha}}{\Gamma[(k+1)\alpha]}
= \sum_{k=0}^{\infty}A^k\frac{t^{(k+1)\alpha-1}}{\Gamma[(k+1)\alpha]}
= t^{\alpha-1}E_{\alpha,\alpha}(At^{\alpha})
\end{equation*}
denote the $\alpha$-exponential matrix function.
\end{definition}

For $\alpha=1$ both functions given in Definitions~\ref{def:ML:2p} and \ref{def:aEM}
are equal and coincide with the classical exponential matrix:
$E_{\alpha}(At^{\alpha})=e_{\alpha}^{At}=\exp(At)$.
In \cite{mt} the following properties are proved:

\begin{proposition}
\label{prop:fromDD}
Let $\alpha>0$. Then,
\begin{enumerate}
\item[(a)] $D^{\alpha}_{0+} e_{\alpha}^{At}=Ae_{\alpha}^{At}$;
\item[(b)] $D^{\alpha}_{T-} \varphi(t)=A\varphi(t)$,
where $\varphi(t)=S(T-t)=e_{\alpha}^{A(T-t)}$;
\item[(c)] $I+\int_0^t Ae_{\alpha}^{A\tau}d\tau=E_{\alpha}(At^{\alpha})$;
\item[(d)] $\frac{d}{dt}E_{\alpha}(At^{\alpha})=Ae_{\alpha}^{A\tau}$, $t>0$;
\item[(e)] $D^{\alpha}_{t_0+}\left(I^{\alpha}_{t_0+}x(t)\right)=x(t)$.
\end{enumerate}
\end{proposition}
Additionally, we need a more general form for the item b)
of Proposition~\ref{prop:fromDD}. In fact this item is also true for
$S(t)=t^{\alpha+\beta-1}E_{\alpha, \alpha+\beta}(At^{\alpha})$.

It is easy to explain the following properties of fractional integrals
of the $\alpha$-exponential function.
If $t_0=0$ we write $I^{\alpha}=I^{\alpha}_{t_0=0+}$.
\begin{itemize}
\item[(a)] Let $k\in\N\cup\{0\}$. Then,
$I^{1-\alpha}\frac{t^{(k+1)\alpha-1}}{\Gamma[(k+1)\alpha]}=\frac{t^{k\alpha}}{\Gamma(k\alpha+1)}$;
\item[(b)] $I^{1-\alpha}e_{\alpha}^{At}=E_{\alpha}(At)$;
\item[(c)] $I^{\beta}e_{\alpha}^{At}=t^{\alpha+\beta-1}E_{\alpha, \alpha+\beta}(At^{\alpha})$;
\item[(d)] $I^{\beta}e^{At}=I^{\beta}e_{\alpha=1}^{At}=t^{\beta}E_{1,\beta+1}(At)$;
\item[(e)] $\left(I^{\beta}e_{\alpha}^{At}\right)_{t=0}=\textbf{0}$, $\beta>1-\alpha$;
\item[(f)] $\left(I^{1-\alpha}e_{\alpha}^{At}a\right)_{t=0}=a$,  $a\in\R^n$.
\end{itemize}

Let $S(t)=t^{\alpha+\beta-1}E_{\alpha, \alpha+\beta}(At^{\alpha})$.
The next result follows from the technical lemmas proved in \cite{mt}.

\begin{lemma}
\label{lemma:intbypart}
Let $0<\alpha<1$
and  $\mu(t)\in I^{\alpha}_{0}(L_q)$, where  $1/p+1/q\leq 1+\alpha$ for $p$ such that
all components of $S(T-t)$ belong to $I^{\alpha}_{T-}(L_p)$.
Then,
\begin{equation*}
\int_0^TS(T-\tau)D^{\alpha}_{0}\mu(\tau)d\tau=\int_0^TAS(T-\tau)\mu(\tau)d\tau \, .
\end{equation*}
\end{lemma}


\section{Systems and the Controllability Gramian}
\label{sec:syst:Gram}

We consider an initialized fractional
linear time-invariant control system
\begin{equation}
\label{eq1}
\tag{$\Sigma_{\alpha,\psi}$}
D^{\alpha}_{0+} x(t)  =  A x(t) + B u(t)\,,
 \ \  t\geq 0\,,
\end{equation}
where $\alpha\in(0,1]$, $x(t)\in \mathbb{R}^n$, $u(t)\in \mathbb{R}^m$,
matrices $A\in \mathbb{R}^{n\times n}$ and $B\in \mathbb{R}^{n\times m}$,
together with an initialization vector $\psi(t)$, $t \leq 0$.
We consider piecewise constant controls $u(\cdot)$.
In \cite{lorenzo} authors explain that the ``state'' of the system is not given
by the dynamic variable vector. This is because of the initialization vector.
The derivative notation is
\[D^{\alpha}_{0+} x(t)=\mbox{}_0d^{\alpha}_tx(t)+\psi(t)\,,\]
where $d^{\alpha}_tx(t)$ is the uninitialized fractional derivative starting
at $t=0$ and $\psi$ is the initialization vector function, determined as
\begin{equation*}
\psi(t)=\lim_{\omega\rightarrow -\infty}\frac{d}{dt}\left(\frac{1}{\Gamma(1-\alpha)}
\int_{\omega}^0\frac{x(\tau)}{(t-\tau)^{\alpha}}d\tau\right)\,.
\end{equation*}
Comparing with \cite{hartley1}, we have $a=-\infty$:
we are involving the history from the left axis.
In the particular case when $x(t)=a$ for $-\infty<t\leq 0$, we obtain the following:
\[\psi(t)=\lim_{\omega\rightarrow -\infty}\frac{d}{dt}\left(\frac{1}{\Gamma(1-\alpha)}
\int_{\omega}^0\frac{a}{(t-\tau)^{\alpha}}d\tau\right)
= -\frac{at^{-\alpha}}{\Gamma(1-\alpha)}\,.
\]

\begin{proposition}
Let $0<\alpha \leq 1$. The forward trajectory of the system $\Sigma_{\alpha,\psi}$
evaluated at $t>0$ has the following form:
\begin{equation*}
\gamma(t,\psi,u)=\int_0^tS(t-\tau)\left(Bu(\tau)-\psi(\tau)\right)d\tau\, ,
\end{equation*}
where  $S(t)= e_{\alpha}^{At}$. Moreover, for $x(t)=a$, $t\in(-\infty,0]$, we have:
\begin{equation}
\label{form}
\gamma(t,a,u)=E_{\alpha}(At^{\alpha})a+\int_0^tS(t-\tau)Bu(\tau)d\tau\, .
\end{equation}
\end{proposition}

\begin{proof}
First part is proved in \cite{hartley1} using the Laplace transform.
We only need to prove here the second part.
We have
\[X(s)=(Is^{\alpha}-A)^{-1}BU(s)-(Is^{\alpha}-A)^{-1}\Psi(s)\, .\]
Applying the inverse Laplace transform for
\[-(Is^{\alpha}-A)^{-1}\Psi(s)=(Is^{\alpha}-A)^{-1}s^{\alpha-1}a\]
we get
\[\mathcal{L}^{-1}\left[-(Is^{\alpha}-A)^{-1}\Psi(s)\right](t)=E_{\alpha}(At^{\alpha})a\,.\]
Hence the formula (\ref{form}) holds.
\end{proof}

We introduce now the more general situation when one has some external measure
of the memory. This is represented by the fractional integral of order $\beta$.
\begin{definition}
Let $ \beta\geq 1-\alpha$. The memory of order $\beta$ of the forward trajectory
$\gamma(t,\psi,u)$ evaluated at time $t> 0$ is defined by
\begin{equation*}
 M_{\beta}(t,\psi,u):=I^{\beta}_{0+}\gamma(t,\psi,u) \, .
\end{equation*}
Additionally, we define the value of the memory of order $\beta$ at $t=0$ as
$M_{\beta}(0,\psi,u):= \lim\limits_{t\rightarrow 0+} I^{\beta}_{0+}\gamma(t,\psi,u)$.
\end{definition}

\begin{proposition}
Let $\beta\geq 1-\alpha$. Then for the system $\Sigma_{\alpha,\psi}$ we have
\begin{equation*}
\label{eq:memory}
M_{\beta}(t,\psi,u)=\int_0^t \Phi_{\beta}(t-\tau)\left(Bu(\tau)-\psi(\tau)\right)d\tau
\  \mbox{for} \ t\geq 0,
\end{equation*}
where
\[\Phi_{\beta}(t)=t^{\alpha+\beta-1}E_{\alpha, \alpha+\beta}(At^{\alpha})
=t^{\alpha+\beta-1}\sum_{k=0}^{\infty}A^k\frac{t^{k\alpha}}{\Gamma((k+1)\alpha+\beta)}\, .\]
Moreover, for $\psi(t)=a$ for $t\in(-\infty,0]$ we have
\begin{equation*}
\label{eq:memory1}
M_{\beta}(t,\psi=a,u)=t^{\beta}E_{\alpha,\beta}(At^{\alpha})a+\int_0^t \Phi_{\beta}(t-\tau)Bu(\tau)d\tau\,.
\end{equation*}
\end{proposition}
\begin{proof}
Using the Laplace transformation we get
\[
\mathcal{L}\left[M_{\beta}(t,\psi,u)\right](s)=\frac{X(s)}{s^{\beta}}
= -\sum_{k=0}^{\infty}A^ks^{-(k+1)\alpha-\beta}\Psi(s)
+\sum_{k=0}^{\infty}s^{-(k+1)\alpha-\beta}A^kBU(s)) \, .
\]
The desired formula follows from the fact that
\begin{equation*}
\begin{split}
\Phi_{\beta}(t) &= \mathcal{L}^{-1}\left[\sum_{k=0}^{\infty}A^ks^{-(k+1)\alpha-\beta}\right](t)\\
&= t^{\alpha+\beta-1}\sum_{k=0}^{\infty}A^k\frac{t^{k\alpha}}{\Gamma((k+1)\alpha+\beta)}\\
&=t^{\alpha+\beta-1}E_{\alpha, \alpha+\beta}(At^{\alpha}) \, .
\end{split}
\end{equation*}
Moreover, for
$\Psi(s)=-s^{\alpha-1}a=\mathcal{L}\left[-\frac{t^{-\alpha}}{\Gamma(1-\alpha)}\right](s)$
we get that
\begin{equation*}
\begin{split}
\mathcal{L}^{-1}&\left[-\sum_{k=0}^{\infty}A^ks^{-(k+1)\alpha-\beta}\Psi(s)\right](t)\\
&= \mathcal{L}^{-1}\left[\sum_{k=0}^{\infty}A^ks^{-(k+1)\alpha-\beta}a\right](t)
=t^{\beta}\sum_{k=0}^{\infty}A^k\frac{t^{k\alpha}}{\Gamma(k\alpha+\beta}a\\
&=t^{\beta}E_{\alpha,\beta}(At^{\alpha})\,.
\end{split}
\end{equation*}
\end{proof}

\begin{definition}
Let $T> 0$. The system $\Sigma_{\alpha,\psi}$ is controllable with memory of order
$\beta\geq 1-\alpha$ on $[0,T]$  if there exists a control
$u$ defined on $[0,T]$ such that \[M_{\beta}(t,\psi,u)_{|t=T}=b\,.\]
\end{definition}

We denote by
\begin{equation*}
Q_T =\int_{0}^{T} (T-t)^{2(1-\alpha-\beta)} \Phi_{\beta}(T-t)BB^* \Phi_{\beta}^*(T-t)dt
\end{equation*}
the \emph{$\beta$-controllability Gramian on the time interval $[0,T]$}
corresponding to the system $\Sigma_{\alpha,\psi}$. As in the classical case,
$Q_{T}$ is symmetric and nonnegative definite.

\begin{theorem}
\label{th1}
Let $T>0$ and $Q_T$ be nonsingular. Then,
\begin{itemize}
\item[(a)] for $b\in \mathbb{R}^n$
the control law
\begin{equation}
\label{control}
\overline{u}(t)=-(T-t)^{2(1-\alpha-\beta)}B^*\Phi_{\beta}^*(T-t)Q_T^{-1}f_T(\psi,b) \, ,
\end{equation}
where $f_T(\psi,b)=-b-\int_0^T \Phi_{\beta}(T-\tau)\psi(\tau)d\tau$,
drives the system  to $b$ in time $T$.
\item[(b)] Among all possible controls from
$L^2_{\alpha}\left([0,T],\R^m\right)$ driving the system to $b$ in time $T$,
the control $\overline{u}$ defined by (\ref{control})
minimizes the modified energy integral
\begin{equation*}
\mathcal{E}^{\alpha,\beta}(u) := \int_0^T \left|(T-t)^{\alpha+\beta-1}u(t)\right|^2 dt \, .
\end{equation*}
\end{itemize}
Moreover,
\[
\mathcal{E}^{\alpha,\beta}(\overline{u}) =
\int_0^T|(T-t)^{\alpha+\beta-1}\overline{u}(t)| ^2dt =<Q_T^{-1}f_T(\psi,b), f_T(\psi,b)>\, ,
\]
where $<\cdot,\cdot>$ denotes the inner product, and
$$
f_T(\psi,b)=-b-\int_0^T \Phi_{\beta}(T-\tau)\psi(\tau)d\tau\, .
$$
\end{theorem}

\begin{proof}
From the form of $\overline{u}$ we directly have that
\begin{equation*}
\begin{split}
M_{\beta}(T,\psi,\overline{u})&=f_T(\psi,b)+b\\
& \quad - \left(\int_0^T (T-t)^{2(1-\alpha-\beta)}
\Phi_{\beta}(T-t)BB^*\Phi_{\beta}^*(T-t)dt\right)Q_T^{-1}f_T(\psi,b) \\
& = f_T(\psi,b)+b-Q_TQ_T^{-1}f_T(\psi,b)\\
& =b \, .
\end{split}
\end{equation*}
Item (a) is proved. Similarly to the classical situation \cite{zabczyk},
to prove item (b) we begin noticing that
\begin{equation*}
\begin{split}
\int_0^T & \left|(T-t)^{\alpha+\beta-1}\overline{u}(t)\right|^2dt \\
&=\int_0^T\left| (T-t)^{1-\alpha-\beta}B^*\Phi_{\beta}^*(T-t)Q_T^{-1}f_T(\psi,b)\right|^2dt\\
&=\int_0^T|T-t|^{2(1-\alpha-\beta)}\left<B^*\Phi_{\beta}^*(t-t)Q_T^{-1}f_T(\psi,b),B^*\Phi_{\beta}^*(T-t)Q_T^{-1}f_T(\psi,b)\right>dt\\
&=\left<\int_0^T |T-t|^{2(1-\alpha-\beta)}\Phi_{\beta}(T-t)BB^*\Phi_{\beta}^*(T-t)dt, Q_T^{-1}f_T(\psi,b)\right>\\
&=\left<Q_TQ_T^{-1}f_T(\psi,b),Q_T^{-1}f_T(\psi,b)\right>\\
&=\left<f_T(\psi,b),Q_T^{-1}f_T(\psi,b)\right>\,.
\end{split}
\end{equation*}
Let us take another control $u$ for which $(T-t)^{\alpha+\beta-1}u(t)$
is square integrable on $[0,T]$ and $M_{\beta}(T,\psi,u)=b$. Then,
\begin{equation*}
\begin{split}
\int_0^T & (T-t)^{2(\alpha+\beta-1)}\left<u(t),\overline{u}(t)\right>dt\\
&=-\int_0^T (T-t)^{2(\alpha+\beta-1)}\left<u(t),(T-t)^{2(1-\alpha-\beta)}B^*\Phi_{\beta}^*(T-t)Q_T^{-1}f_T(\psi,b)\right>dt\\
&=-\int_0^T\left<u(t),B^*\Phi_{\beta}^*(T-t)Q_T^{-1}f_T(\psi,b)\right>dt \\
&= \left<f_T(\psi,b),Q_T^{-1}f_T(\psi,b)\right>\,.
\end{split}
\end{equation*}
Hence, \[\int_0^T(T-t)^{2(1-\alpha-\beta)}\left<u(t),\overline{u}(t)\right>dt
=\int_0^T(T-t)^{2(1-\alpha-\beta)}\left<\overline{u}(t),\overline{u}(t)\right>dt\]
and we obtain
\begin{multline*}
\int_0^T(T-t)^{2(1-\alpha-\beta)}|u(t)|^2dt\\
=\int_0^T(T-t)^{2(1-\alpha-\beta)}|\overline{u}(t)|^2dt+\int_0^T(T-t)^{2(1-\alpha-\beta)}|u(t)-\overline{u}(t)|^2 dt \, ,
\end{multline*}
which gives the intended minimality property for the integral.
 \end{proof}

\begin{example}
Let $\Sigma_{\alpha=0.5,\psi}$ be the following system in $\R^2$:
\[
\begin{cases}
D^{0.5}_{0+} x_1(t)=x_2(t) \, , \\
D^{0.5}_{0+} x_2(t)=u(t)\,.
\end{cases}
\]
Let us take
$b=\left(\begin{array}{c} 0 \\  0 \end{array}\right)$.
Since $A=\left(\begin{array}{cc} 0 & 1\\ 0 & 0\end{array}\right)$,
$B=\left(\begin{array}{c} 0 \\  1 \end{array}\right)$,
and $A^2=\textbf{0}$, we obtain for $t_0=0$
the formula for the solution with the initialization $\psi$
 in the following form:
\[\gamma(t,\psi,u)=\int_0^t \left(\begin{array}{c} 1 \\
\frac{1}{\sqrt{\pi(t-\tau)}} \end{array}\right)\left(u(\tau)-\psi(\tau)\right)d\tau\, . \]
Let $\beta\geq 0.5$. Then,
\[\Phi_{\beta}(t)=\left(\begin{array}{cc} \frac{t^{\beta-0.5}}{\Gamma(0.5+\beta)}
& \frac{t^{\beta}}{\Gamma(1+\beta)} \\ 0 &  \frac{t^{\beta-0.5}}{\Gamma(0.5+\beta)}
\end{array}\right)\,.\]
For simplicity let us take $\beta=0.5$. Then,
\[\Phi_{\beta=0.5}(t)=\left(\begin{array}{cc} 1 &
\frac{2\sqrt{t}}{\sqrt{\pi}} \\ 0 &  1\end{array}\right)\]
and  the Gramian is the following:
\[Q_{T,\beta=0.5}=\left(\begin{array}{cc}\frac{2T^2}{\pi} & \frac{4T^{1.5}}{3\sqrt{\pi}} \\
\frac{4T^{1.5}}{3\sqrt{\pi}} & T\end{array}\right)\,.\]
Thus,
$\overline{u}(t)=\left[\frac{\sqrt{T-t}}{\Gamma(3/2)},1\right]Q_T^{-1}f_T(\psi,b)$.
Let us take the end point $b$ equal to the zero vector and $\psi_2(t)\equiv 0$. Then,
$\overline{u}(t)=\left(\frac{6\sqrt{\pi}}{T^{3/2}
-\frac{9\sqrt{\pi}\sqrt{T-t}}{T^2}}\right)\int_0^T\psi_1(\tau)d\tau$,
and the norm \[\varepsilon^{0.5,0.5}(\overline{u})=
\left(\int_0^T\psi_1(\tau)d\tau\right)^2\int_0^T\left(\frac{6\sqrt{\pi}}{T^{3/2}}
-\frac{9\sqrt{\pi}\sqrt{T-t}}{T^2}\right)dt\]
with $\overline{u}(t)=\frac{9\pi}{2T^2}\left(\int_0^T\psi_1(\tau)d\tau\right)^2$.
\end{example}


\section{Steering Laws}
\label{sec:rank}

In the classical theory it is quite easy to explain
why the rank condition of matrix $B$
equal to $n$ is sufficient for controllability.
In the case of fractional order systems
we need to use a special function.

\begin{proposition}
\label{prop:control}
Let $\mbox{rank}\, B=n$, matrix $B^+$ be such that
$BB^+=I$, and $g_{\beta}(\cdot)$ be the matrix function
such that $\Phi_{\beta}(t)g_{\beta}(t)=I$, $t>0$,
and $\lim\limits_{t\rightarrow 0+}\Phi_{\beta}(t)g_{\beta}(t)=I$.
Then the control
$$
\widehat{u}(t)=\frac{\Gamma(\alpha+\beta)}{T}B^+g(T-t)\left(b+\int_0^T \Phi_{\beta}(T-\tau)\psi(\tau)d\tau\right)\, ,
\quad t\in[0,T] \, ,
$$
transfers system $\Sigma_{\alpha,\psi}$ to $b$ in time $T> 0$.
\end{proposition}

\begin{proof}
It follows by direct calculation:
\[\gamma(T,\psi,\widehat{u})=-\int_0^T \Phi_{\beta}(T-\tau)\psi(\tau)d\tau\]
\[+\frac{1}{T}\int_0^T\Phi_{\beta}(T-t)BB^+
g_{\beta}(T-t)\left(b+\int_0^T \Phi_{\beta}(T-\tau)\psi(\tau)d\tau\right)dt=b\,.\]
\end{proof}

\begin{example}
Let us consider the system
\[\Sigma_{\alpha,\psi}: \quad D_{0+}^{\alpha}x(t)=u(t) \, , \]
where $D_{0+}^{\alpha}x(t)=d^{\alpha}_{0+}x(t)+\psi(t)$.
Let $b\in\R$, $T> 0$. Then,
\begin{equation*}
\begin{gathered}
\Phi_{\beta}(t) = \frac{t^{\alpha+\beta-1}}{\Gamma(\alpha+\beta)} \, , \quad B^+ = B=1 \, , \\
M(T,\psi,u) =\frac{1}{\Gamma(\alpha+\beta)}\int_0^T(T-t)^{\alpha+\beta-1}\left(u(t)-\psi(t)\right)dt\, .
\end{gathered}
\end{equation*}
The steering law
\[\widehat{u}(t)=\frac{\Gamma(\alpha+\beta)}{T}(T-t)^{1-\alpha-\beta}\left(b+\int_0^T \Phi_{\beta}(T-\tau)\psi(\tau)d\tau\right)\]
transfers the system to $b$.
\end{example}

If the rank condition is satisfied, the control $\overline{u}$
given by (\ref{control}) drives the system to $b$ at time $T$.
Our goal now is to find another formula for the steering control by using the matrix $\left[A|B\right]$
instead of the controllability matrix $Q_T$. It is a classical result that if $\mbox{rank}\,\left[A|B\right]=n$
then there exists a matrix $K\in M(mn,n)$ such that $\left[A|B\right]K=I\in M(n,n)$ or, equivalently,
there are matrices $K_1$, $K_2, \ldots K_n \in M(m,n)$ such that
\begin{equation*}
BK_1+ABK_2+\cdots +A^{n-1}BK_n=I\,.
\end{equation*}
In \cite{mt} we introduced a convenient notation for composition
of fractional derivatives with the same order $\alpha$,
where $\alpha\in(0,1]$. Let us recall here the notation
for fractional Riemann-Liouville derivatives. Let
$R^{\alpha,0}_{0+}\mu(t)=\mu(t)$.
We define, recursively,
$R^{\alpha,j+1}_{0+}\mu(t):=D^{\alpha}_{0+}\left(R^{\alpha,j}_{0+}\mu(t)\right)$,
$j\in\N$.

\begin{theorem}
\label{thm:th2}
Let $\beta\geq 0$ and $\mbox{rank}\,\left[A|B\right]
=\mbox{rank}\,\left[B,AB, \ldots, A^{n-1}B\right]=n$ for system $\Sigma_{\alpha,\psi}$.
Let $p$ be such that $\Phi_{\beta}(T-t)\in I^{\alpha}_{T-}(L_p)$
and let $\varphi$ be a real function given on $[0,T]$ such that
\begin{enumerate}
\item[(1)] $\int_0^T \varphi(t)dt=1$;
\item[(2)] $R^{\alpha,j}_{0+}\mu(t)\in I^{\alpha}_{0+}(L_q)$
for $j=0,\ldots, n-1$, where
\[\mu(t)=g_{\beta}(t)\left(b+\int_0^T \Phi_{\beta}(T-\tau)\psi(\tau)d\tau\right)\varphi(t)\, ,\]
and $\Phi_{\beta}(T-t)g_{\beta}(t)=I_n$ for $t\in [0,T]$, while $1/p+1/q\leq 1+\alpha$.
\end{enumerate}
Then, the control
\begin{equation*}
\hat{u}(t)=K_1\mu(t)+K_2D^{\alpha}_{0+}\mu(t)+\cdots+K_nR^{\alpha, n-1}_{0+}\mu(t), \ \ t\in[0,T] \, ,
\end{equation*}
transfers the system $\Sigma_{\alpha, \psi}$ to the value $b$ at time $T>0$.
\end{theorem}

\begin{proof}
Using $j-1$ times the formula (\ref{eq:intbypart}) of integration by parts
and Lemma~\ref{lemma:intbypart}, we get for $j=1,\ldots, n$ that
\begin{equation*}
\int_0^T \Phi_{\beta}(T-t)BK_jR^{\alpha, j-1}_{0+}\mu(t)dt
=\int_0^T \Phi_{\beta}(T-t)A^{j-1}BK_j\mu(t)dt\,.
\end{equation*}
Then,
\[\int_0^T \Phi_{\beta}(T-t)B\hat{u}(t)dt=\int_0^T \Phi_{\beta}(T-t)B\mu(t)dt\]
and finally
\[\gamma(T,\psi,\hat{u})= -\int_0^T \Phi_{\beta}(T-\tau)\psi(\tau)d\tau\]
\[+\int_0^T
\Phi_{\beta}g_{\beta}(t)\left(b+\int_0^T \Phi_{\beta}(T-\tau)\psi(\tau)d\tau\right)\varphi(t)dt=b\,.\]
\end{proof}


\section{Conclusions}
\label{sec:cnc}

In this work we develop the theoretical aspects
regarding optimality and steering laws for fractional systems
$\Sigma_{\alpha,\psi}$ with Riemman-Liouville derivatives of order $\alpha$
and a given initialization function $\psi$.
The novelty of the theory here introduced
is the ``memory state point of view'' for
studying fractional systems. Main results provide
explicit steering laws that drive the system from
the initial memory value characterized by $\psi$ to a desired
memory state of the system. Moreover, we obtained
the optimal law in the sense of minimization of the modified energy.

The present work opens several possibilities of future research.
Interesting open questions include the study of memory observability conditions
for initialized systems and realization issues. The realization
question is particularly interesting. Roughly speaking,
after solving a fractional system of order $\alpha$
we neglect that we know the true state,
hiding it under an additional integral of a different fractional order.
From the point of view of transfer functions, it means that we do not know of what type,
in the sense of the fractional order and class of derivative,
a realization can be done. One possible approach is to realize
the transfer function using not only the direct dynamic
but coming back to the dynamic of the additional integral of the state.
This research is in its beginning, and it will be developed in the future.


\section*{Acknowledgements}

DM was supported by BUT grant
S/WI/1/08; DT by the R\&D unit CIDMA,
via FCT and the EC fund FEDER/POCI 2010.
The authors would like to express their gratitude to two anonymous
referees, for several relevant and stimulating remarks contributing
to improve the quality of the paper, in particular
for references regarding the initialization issue.



\end{document}